\documentclass[reqno,a4paper]{amsart}
\usepackage{amsmath}
\usepackage{amssymb,amsfonts}
\usepackage{amsthm}
\usepackage{enumerate}
\usepackage[mathscr]{eucal}
\usepackage{eqlist}

\theoremstyle{plain}
\newtheorem{theorem}{Theorem}[section]

\newtheorem{lemma}{Lemma}[section]

\theoremstyle{definition}

\newtheorem{remark}{\textnormal{\textbf{Remark}}}

\theoremstyle{remark}

\numberwithin{equation}{section}

\begin{document}

\title[On the power values of the sum of three square in arith. prog.]{On the power values of the sum of three squares in arithmetic progression}

\author[M. Le and G. Soydan]{Maohua Le and G\"{o}khan Soydan}

\address{{\bf Maohua Le}\\
	Institute of Mathematics, Lingnan Normal College\\
	Zhangjiang, Guangdong, 524048 China}

\email{lemaohua2008@163.com}

\address{{\bf G\"{o}khan Soydan} \\
	Department of Mathematics \\
	Bursa Uluda\u{g} University\\
	G\"{o}r\"{u}kle Campus, 16059 Bursa, Turkey}
\email{gsoydan@uludag.edu.tr }

\newcommand{\acr}{\newline\indent}

\thanks{}

\subjclass[2010]{11D41,11J86}
\keywords{polynomial Diophantine equation, power sums, primitive divisors of Lehmer sequences, Baker's method}

\begin{abstract}
In this paper, using a deep result on the existence of primitive divisors of Lehmer numbers due to Y. Bilu, G. Hanrot and P. M. Voutier, we first give an explicit formula for all positive integer solutions of the Diophantine equation $(x-d)^2+x^2+(x+d)^2=y^n$ (*) when $n$ is an odd prime and $d=p^r$, $p>3$ a prime. So this improves the results on the papers of A. Koutsianas and V. Patel \cite{KP} and A. Koutsianas \cite{Kou}. Secondly, under the assumption of our first result, we prove that (*) has at most one solution $(x,y)$. Next, for a general $d$, we prove the following two results: (i) if every odd prime divisor $q$ of $d$ satisfies $q\not\equiv \pm 1 \pmod{2n},$ then (*) has only the solution $(x,y,d,n)=(21,11,2,3)$. (ii) if $n>228000$ and $d>8\sqrt{2}$, then all solutions $(x,y)$ of (*) satisfy $y^n<2^{3/2}d^3$.
\end{abstract}

\maketitle

\section{Introduction}\label{sec:1}
Let $\mathbb{Z}$, $\mathbb{N}$ and $\mathbb{Q}$ be the sets of all integers, positive integers and rational numbers respectively. Let $k,n$ be fixed positive integers. The study of the polynomial Diophantine equation in the form of
\begin{equation}\label{eq.1.1}
1^{k}+2^{k}+...+x^{k}=y^n, \quad x,y\in\mathbb{N},\quad n\geq 2
\end{equation}
has been going on for more than a hundred years. In 1875, the classical question of E. Lucas \cite{Lu} was whether equation \eqref{eq.1.1} has only the solutions $x=y=1$ and $x=24$, $y=70$ for $(k,n)=(2,2)$. In 1918, G. N. Watson \cite{W} solved equation \eqref{eq.1.1} with $(k,n)=(2,2)$. In 1956,  J. J. Sch\"{a}ffer \cite{Sch} considered  equation \eqref{eq.1.1}. He showed, for $k\geq 1$ and $n\geq 2$, that \eqref{eq.1.1} possesses at most finitely many solutions in positive integers $x$ and $y$, unless
\begin{equation}  \label{eq.1.2}
(k,n) \in \{(1,2),(3,2),(3,4),(5,2)\},
\end{equation}
where, in each case, there are infinitely many such solutions. J. J. Sch\"{a}ffer's conjectured that \eqref{eq.1.1} has the unique non-trivial (i.e. $(x,y)\neq (1,1)$) solution, namely $(k,n,x,y)=(2,2,24,70)$. The correctness of this conjecture has been proved for some cases (see, e.g., \cite{BGP}, \cite{BHMP}, \cite{GP}, \cite{H}, \cite{JPW}, \cite{Pi1}, \cite{Pi2}). But, it has not been proved completely yet.

A more general case is to consider the Diophantine equation
\begin{equation}  \label{eq.1.3}
	(x+1)^k+(x+2)^k+ ... + (x+r)^k=y^{n} \quad x,y\in\mathbb{Z},\quad k,n\geq 2.
\end{equation}

In 2013, Z. Zhang and M. Bai \cite{ZB} solved the equation \eqref{eq.1.3} with $k=2$ and $r=x$. In 2014, the equation
\begin{equation}\label{eq.1.4}
(x-1)^k+x^k+(x+1)^k=y^n  \quad x, y \in\mathbb{Z},\quad n \geq 2,
\end{equation}
was solved completely by Z. Zhang \cite{Zh1} for $k=2, 3, 4$ (Actually, firstly, J. W. S. Cassels \cite{Cas} considered the equation \eqref{eq.1.4} in 1985, and he proved that $x=0,1,2,24$ are the only integer solutions to this equation for $k=3$ and $n=2$) and in 2016, M. A. Bennett, V. Patel and S. Siksek \cite{BPS}  extended Z. Zhang's result, completely solving equation \eqref{eq.1.4} in the cases $k=5$ and $k= 6$. The same year, M. A. Bennett, V. Patel and S. Siksek \cite{BPS2} considered the equation \eqref{eq.1.3}. They gave the integral solutions to the equation \eqref{eq.1.3} using linear forms in logarithms, sieving and Frey curves where $k=3$, $2\leq r \leq 50,$ $x\geq1$ and $n$ is prime.

Let $k\geq2$ be even, and let $r$ be a fixed non-zero integer. In 2017, V. Patel and S. Siksek \cite{PS} showed that
for almost all $d\geq2$ (in the sense of natural density), the equation
\begin{equation*}
x^k +(x+r)^k +...+(x+(d-1)r)^k =y^n, \quad x, y \in\mathbb{Z},  \quad n \geq 2
\end{equation*}
has no solutions. Let $\ell \geq2$ be a fixed integer such that $\ell$ even. The same year, the second author \cite{S} considered the equation
\begin{equation}\label{eq.1.5}
(x+1)^k+(x+2)^k+...+(\ell x)^k =y^n, \quad x, y\in\mathbb{Z}  \quad n \geq 2.
\end{equation}
He proved that the equation \eqref{eq.1.5} has only finitely many solutions where $x,y\ge 1$, $k\neq 1,3$. He also showed that the equation \eqref{eq.1.5} has infinitely many solutions with $n\geq 2$ and $k=1,3$. In 2018, A. B\'{e}rczes, I. Pink, G. Sava\c{s} and the second author \cite{BPSS} considered the equation \eqref{eq.1.5} with $\ell=2$. They proved that the equation \eqref{eq.1.5} has no solutions where $2\leq x\leq 13$, $k\geq 1$, $\ell=2$, $y\geq 2$ and $n\geq 3$. Recently, D. Bartoli and the second author \cite{BarSoy} proved that all the solutions of the equation \eqref{eq.1.5} with $x,y\geq1, n\geq2, k\neq3$ and $\ell$ odd satisfy $\max\{x,y,n\}<C$ where $C$ is an effectively computable constant depending only on $k$ and $\ell$. So, the remaining case for the equation \eqref{eq.1.5} was covered by them.

Finding perfect powers that are sums of terms in an arithmetic progression has received much interest; recent contributions can be also found in \cite{G}, \cite{GP2}, \cite{BK}, \cite{DDKT}.

Now we consider a generalization of equation \eqref{eq.1.4}. Let $d$ be fixed positive integer. In 2017-2019, Z. Zhang \cite{Zh2} A. Koutsianas and V. Patel \cite{KP} studied the integer solutions of the following equation
\begin{equation}\label{eq.1.6}
(x-d)^{k}+x^{k}+(x+d)^{k}=y^n,\quad x,y\in \mathbb{Z},\quad n\geq2
\end{equation}
for the cases $k=4$ and $k=2$, respectively.
Z. Zhang gave some results on the equation \eqref{eq.1.6} with $k=4$ by using modular approach. A. Koutsianas and V. Patel \cite{KP} gave  all non-trivial primitive solutions to equation \eqref{eq.1.6} where $k=2$, $n$ is prime and $d\le 10^4$. (According to the terminology of \cite{KP}, an integer solution $(x,y)$ of \eqref{eq.1.6} is said to be primitive if $\gcd(x,y)=1$. This is equivalent to $x,y,d$ being pairwise coprime. A solution where $xy=0$ is called a trivial solution). They used the characterization of primitive divisors in Lehmer sequences due to Y. F. Bilu, G. Hanrot and P. M. Voutier \cite{BHV}, then A. A. Garcia and V. Patel \cite{GaP} showed that the only solutions to the equation \eqref{eq.1.6} with $n\ge5$ a prime, $k=3$, $\gcd(x,d)=1$ and $0<d \le 10^6$ are the trivial ones satifying $xy=0$.

Recently, A. Koutsianas \cite{Kou} studied the equation \eqref{eq.1.6} with $k=2$ for an infinitely family of $d$ which is an extension of \cite{KP}. In \cite{Kou}, all solutions $(x,y)$ of the Diophantine equation
\begin{equation}\label{maineq}
	(x-d)^{2}+x^{2}+(x+d)^{2}=y^n,\quad x,y\in \mathbb{N},\quad n\ge 2,\quad \gcd(x,y)=1,
\end{equation}
are given with the following table where $d=p^r$ with $r\ge 0$, $p$ a prime and $p\le 10^4$. 

\begin{table}[ht]
	\caption{}
	\begin{center}
		\begin{tabular}{|l|l|}
			
			\hline
			$p$ & $(x,y,r,n)$ \\ \hline
			$2$ & $(21,11,1,3)$ \\ \hline
			$7$ & $(3,5,1,3)$ \\ \hline
			$79$ & $(63,29,1,3)$ \\ \hline
			$223$ & $(345,77,1,3)$ \\ \hline
			$439$ & $(987,149,1,3)$ \\ \hline
			$727$ & $(2133,245,1,3)$ \\ \hline
			$1087$ & $(3927,365,1,3)$ \\ \hline
			$3109$ & $(627,29,1,5)$ \\ \hline
			$3967$ & $(27657,1325,1,3)$ \\ \hline
			$4759$ & $(36363,1589,1,3)$ \\ \hline
			$5623$ & $(46725,1877,1,3)$ \\ \hline
			$8647$ & $(89187,2885,1,3)$ \\ \hline		
		\end{tabular}
	\end{center}
\end{table}

However, the Table 1 at least omits the solution $(x,y,d,r,n)=(13,5,197,1,7)$ of \eqref{maineq} with $p\le 10^4$. 

\bigskip

In this paper, extending the results in \cite{Kou} and \cite{KP}, we first consider the Diophantine equation \eqref{maineq} where
\begin{equation}\label{eq.xx}
d=p^r \,\,\,\text{with}\,\,\, r\in \mathbb{N}.
\end{equation}
We prove the following two results:  

\begin{theorem}\label{theo:1.1}
Let $n$ be an odd prime, and let $d$ be satisfied as in \eqref{eq.xx}. If $(x,y)$ is a solution of \eqref{maineq}, then $p>3$ and there exists a constant $X_1\in\mathbb{N}$ such that
\begin{equation}\label{eq.1.7}
d=\left\lvert \sum\limits_{i=0}^{(n-1)/2}\binom {n} {2i+1}(3X_1^2)^{(n-1)/2-i}(-2)^i\right\rvert.
\end{equation}
Moreover, if \eqref{eq.1.7} holds, then the solution $(x,y)$ can be expressed as
\begin{equation}\label{eq.1.8}
x=X_1\left\lvert \sum\limits_{i=0}^{(n-1)/2}\binom {n} {2i}(3X_1^2)^{(n-1)/2-i}(-2)^i\right\rvert,\,\,y=3X_1^2+2.
\end{equation}
\end{theorem}
\begin{remark}
Theorem \ref{theo:1.1} gives the missing solution $(x,y,d,r,n)=(13,5,197,1,7)$ in \cite{Kou} where $X_1=1$ and $n=7$.
\end{remark}
\begin{theorem}\label{theo:1.2}
Under assumption of Theorem \ref{theo:1.1}, \eqref{maineq} has at most one solution $(x,y)$.
\end{theorem}
Please note that in \cite{Kou}, while all solutions $(x,y)$ of \eqref{maineq} are given where $d=p^r$ with $r\ge 0$, $p$ a prime and $p\le 10^4$, Theorem \ref{theo:1.1} gives an explicit formula to find all solutions $(x,y)$ of \eqref{maineq} for all $d=p^r$ with $r\in\mathbb{N}$.

Next, for a general $d$, we prove the following two results:
\begin{theorem}\label{theo:1.3}
If $n$ is an odd prime and every odd prime divisor $q$ of $d$ satisfies $q\not\equiv \pm 1 \pmod{2n},$ then \eqref{maineq} has only the solution $(x,y,d,n)=(21,11,2,3)$.
\end{theorem}

\begin{theorem}\label{theo:1.4}
If $n>228000$ and $d>8\sqrt{2}$, then all solutions $(x,y)$ of \eqref{maineq} satisfy $y^n<2^{3/2}d^3$.	
\end{theorem}

\section{Proof of Theorem \ref{theo:1.1}}
Let $D_1,D_2,k$ be fixed positive integers such that $\min\{D_1,D_2\}>1$, $2\nmid k$ and $\gcd(D_1,D_2)=\gcd(D_1D_2,k)=1$, and let $h(-4D_1D_2)$ denote the class number of positive binary quadratic primitive forms with discriminant $-4D_1D_2$.

\begin{lemma}\label{lem.2.1}
If the equation
\begin{equation*}
D_1X^2+D_2Y^2=k^Z,\,\,X,Y,Z\in\mathbb{Z},\,\,\gcd(X,Y)=1,\,Z>0
\end{equation*}
has solutions $(X,Y,Z)$, then its every solution $(X,Y,Z)$ can be expressed as
\begin{equation*}
\begin{aligned}
&Z=Z_1t, \, t\in\mathbb{N},\,\,2\nmid t,\\
&X\sqrt{D_1}+Y\sqrt{-D_2}=\lambda_1(X_1\sqrt{D_1}+\lambda_2Y_1\sqrt{-D_2})^t, \, \lambda_1, \lambda_2 \in\{1,-1\}.
\end{aligned}
\end{equation*}
where $X_1,Y_1,Z_1$ are positive integers such that
\begin{equation*}
D_1X_1^2+D_2Y_1^2=k^{Z_1},\,\,\gcd(X_1,Y_1)=1
\end{equation*}
and $h(-4D_1D_2)\equiv 0 \pmod{2Z_1}$.
\end{lemma}
\begin{proof}
This is special case of Theorems 1 and 3 of \cite{Le} for $D<0$ and $D_1>1$.
\end{proof}
\begin{lemma}\label{lem.2.2}
If \eqref{maineq} has solutions $(x,y)$, then $2\nmid n$ and its every solution $(x,y)$ can be expressed as 
\begin{equation}\label{eq.2.1}
x\sqrt{3}+d\sqrt{-2}=\lambda_{1}(X_{1}\sqrt{3}+\lambda_{2}Y_{1}\sqrt{-2})^{n}, \quad \lambda_{1},\lambda_{2}\in\{\pm 1\},
\end{equation}
\begin{equation}\label{eq.2.2}
y=3X_1^2+2Y_1^2,\,X_1,Y_1\in\mathbb{N},\,\,\gcd(X_1,Y_1)=1.
\end{equation}
\end{lemma}
\begin{proof}
We now assume that $(x,y)$ is a solution of \eqref{maineq}. Then we have
\begin{equation}\label{eq.2.3}
3x^2+2d^2=y^n.
\end{equation}
Since $n>2$ and $\gcd(x,y)=1$, by \eqref{eq.2.3}, we get
\begin{equation}\label{eq.2.4}
2\nmid x,\,\,2\nmid y,\,\,3\nmid y,\,\,\gcd(x,d)=1. 
\end{equation}
Hence, we see from \eqref{eq.2.3} and \eqref{eq.2.4} that $\gcd(6,y)=1$ and the equation
\begin{equation}\label{eq.2.5}
3X^2+2Y^2=y^Z,\,\,X,Y,Z \in\mathbb{Z},\,\,\gcd(X,Y)=1,\,\,Z>0,
\end{equation}
has a solution
\begin{equation}\label{eq.2.6}
(X,Y,Z)=(x,d,n).
\end{equation}
Applying Lemma \ref{lem.2.1} to \eqref{eq.2.5} and \eqref{eq.2.6}, we have
\begin{equation}\label{eq.2.7}
n=Z_1t, \,\, t\in\mathbb{N},\,\,2\nmid t,
\end{equation}
\begin{equation}\label{eq.2.8}
x\sqrt{3}+d\sqrt{-2}=\lambda_1(X_1\sqrt{3}+\lambda_2Y_1\sqrt{-2})^t, \, \, \lambda_1, \lambda_2 \in\{1,-1\},	
\end{equation}
where $X_1,Y_1,Z_1$ are positive integers such that
\begin{equation}\label{eq.2.9}
3X_1^2+2Y_1^2=y^{Z_1},\,\, \gcd(X_1,Y_1)=1,
\end{equation}
and
\begin{equation}\label{eq.2.10}
h(-24)\equiv 0 \pmod{2Z_1}.
\end{equation}
Further, since $h(-24)=2$, by \eqref{eq.2.10}, we get $Z_1=1.$ Hence, by \eqref{eq.2.7}, we have $t=n$ and $2\nmid n.$ Furthermore, by \eqref{eq.2.8} and \eqref{eq.2.9}, we obtain \eqref{eq.2.1} and \eqref{eq.2.2} respectively. Thus, Lemma is proved.
\end{proof}

Let $\alpha$, $\beta$ be algebraic integers. If $(\alpha+\beta)^{2}$ and $\alpha\beta$ are nonzero coprime integers and $\alpha/\beta$ is not a root of unity, then $(\alpha,\beta)$ is called a \textit{Lehmer pair}. Further, let $A=(\alpha+\beta)^{2}$ and $C=\alpha\beta$. Then we have
\begin{equation*}
\alpha=\frac{1}{2}(\sqrt{A}+\lambda\sqrt{B}), \quad \beta=\frac{1}{2}(\sqrt{A}-\lambda\sqrt{B}), \quad \lambda\in\{\pm1\},
\end{equation*}
where $B=A-4C$. Such $(A,B)$ is called the parameters of Lehmer pair $(\alpha,\beta)$. Two Lehmer pairs $(\alpha_{1},\beta_{1})$ and $(\alpha_{2},\beta_{2})$ are called equivalent if $\alpha_{1}/\alpha_{2}=\beta_{1}/\beta_{2}\in\{\pm 1,\pm \sqrt{-1}\}$. Obviously, if $(\alpha_{1},\beta_{1})$ and $(\alpha_{2},\beta_{2})$ are equivalent Lehmer pairs with parameters $(A_{1},B_{1})$ and $(A_{2},B_{2})$ respectively, then $(A_{2},B_{2})=(\varepsilon A_{1},\varepsilon B_{1})$, where $\varepsilon\in\{\pm1\}$. For a fixed Lehmer pair $(\alpha,\beta)$, one defines the corresponding \textit{sequence of Lehmer numbers} by
\begin{equation}\label{eq.2.11}
	L_{m}(\alpha,\beta)=\left\{
	\begin{array}{ll}
		\dfrac{\alpha^{m}-\beta^{m}}{\alpha-\beta}, & \textrm{ if $2\nmid m$,} \\
		\dfrac{\alpha^{m}-\beta^{m}}{\alpha^{2}-\beta^{2}}, & \textrm{ if $2|m$, $m\in\mathbb{N}$.} \\
	\end{array}
	\right.
\end{equation}
Then, Lehmer numbers $L_{m}(\alpha,\beta)$ $(m=1,2,...)$ are nonzero integers. Further, for equivalent Lehmer pairs $(\alpha_{1},\beta_{1})$ and $(\alpha_{2},\beta_{2})$, we have $L_{m}(\alpha_{1},\beta_{1})=\pm L_{m}(\alpha_{2},\beta_{2})$ for any $m$. A prime $q$ is called a \textit{primitive divisor} of the Lehmer number $L_{m}(\alpha,\beta)$ $(m>1)$, if $q|L_{m}(\alpha,\beta)$ and $q\nmid ABL_{1}(\alpha,\beta)\cdots L_{m-1}(\alpha,\beta)$, where $(A,B)$ is the parameters of Lehmer pair $(\alpha,\beta)$. For a fixed positive integer $m$, a Lehmer pair $(\alpha,\beta)$ such that $L_{m}(\alpha,\beta)$ has no primitive divisor will be called \textit{$m$-defective} Lehmer pair. Further, a positive integer $m$ is called totally \textit{non-defective} if no Lehmer pair is $m$-defective.

\begin{lemma}[\cite{V}]\label{lem.2.3}
Let $m$ be such that $6<m\leq 30$ and $m\neq 8,10,12$. Then up to equivalence, all parameters $(A,B)$ $(A>0)$ of $m$-defective Lehmer pairs are given as follows:\\
	
$(i)$ $m=7$, $(A,B)=(1,-7), (1,-19), (3,-5), (5,-7), (13,-3), (14,-22).$\\
	
$(ii)$ $m=9$, $(A,B)=(5,-3), (7,-1), (7,-5).$\\
	
$(iii)$ $m=13$, $(A,B)=(1,-7).$\\
	
$(iv)$ $m=14$, $(A,B)=(3,-13), (5,-3), (7,-1), (7,-5), (19,-1), (22,-14).$\\
	
$(v)$ $m=15$, $(A,B)=(7,-1), (10,-2).$\\
	
$(vi)$ $m=18$, $(A,B)=(1,-7), (3,-5), (5,-7).$\\
	
$(vii)$ $m=24$, $(A,B)=(3,-5), (5,-3).$\\
	
$(viii)$ $m=26$, $(A,B)=(7,-1).$\\
	
$(ix)$ $m=30$, $(A,B)=(1,-7), (2,-10).$\\
\end{lemma}

\begin{lemma}[\cite{BHV}]\label{lem.2.4}
Every positive integer $m$ with $m>30$ is totally non-defective.
\end{lemma}
\textbf{Proof of Theorem \ref{theo:1.1}} We now assume that $(x,y)$ is a solution of \eqref{maineq}. Then, $x,y$ and $d$ satisfy \eqref{eq.2.3}. If $p=3$, then from \eqref{eq.xx} and \eqref{eq.2.3} we get $3\mid y$, which contradicts \eqref{eq.2.4}. So we have $p>3$.

By Lemma \ref{lem.2.2}, there exist positive integers $X_1$ and $Y_1$ satisfying \eqref{eq.2.1} and \eqref{eq.2.2}. By \eqref{eq.2.1}, we have
\begin{equation}\label{eq.2.12}
x=X_1\left\lvert \sum\limits_{i=0}^{(n-1)/2}\binom {n} {2i}(3X_1^2)^{(n-1)/2-i}(-2Y_1^2)^i\right\rvert,
\end{equation}
and
\begin{equation}\label{eq.2.13}
d=Y_1\left\lvert \sum\limits_{i=0}^{(n-1)/2}\binom {n} {2i+1}(3X_1^2)^{(n-1)/2-i}(-2Y_1^2)^i\right\rvert.
\end{equation}
Since $d$ satisfies \eqref{eq.xx}, by \eqref{eq.2.13}, we get
\begin{equation}\label{eq.2.14}
Y_1=p^s,\,\,s\in\mathbb{Z},\,\,0\le s\le r,
\end{equation}
and
\begin{equation}\label{eq.2.15}
\left\lvert \sum\limits_{i=0}^{(n-1)/2}\binom {n} {2i+1}(3X_1^2)^{(n-1)/2-i}(-2Y_1^2)^i\right\rvert=p^{r-s}.
\end{equation}

Let
\begin{equation}\label{eq.2.16}
\alpha=X_1\sqrt{3}+Y_1\sqrt{-2},\,\,\beta=X_1\sqrt{3}-Y_1\sqrt{-2}.
\end{equation}
By \eqref{eq.2.2} and \eqref{eq.2.16}, we have
\begin{equation}\label{eq.2.17}
\alpha+\beta=2X_1\sqrt{3},\,\,\alpha-\beta=2Y_1\sqrt{-2},\,\,\alpha\beta=y.
\end{equation}
Notice that $y\ge 5$ by \eqref{eq.2.2}, and $\alpha/\beta$ satisfies
\begin{equation}\label{eq.2.17-A}
y\left(\dfrac{\alpha}{\beta}\right)^2-2(3X_1^2-2Y_1^2)\dfrac{\alpha}{\beta}+y=0
\end{equation}
with $\gcd(y,2(3X_1^2-2Y_1^2))=1$. This implies that $\alpha/\beta$ is not a root of unity. Hence, we see from \eqref{eq.2.4}, \eqref{eq.2.16} and \eqref{eq.2.17} that $(\alpha,\beta)$ is a Lehmer pair with the parameters
\begin{equation}\label{eq.2.18}
(A,B)=(12X_1^2,-8Y_1^2).
\end{equation}
Further, let $L_m(\alpha,\beta)$ $(m=1,2,\cdots)$ be the corresponding Lehmer numbers. By \eqref{eq.2.11} and \eqref{eq.2.16}, we have
\begin{equation}\label{eq.2.19}
\sum\limits_{i=0}^{(n-1)/2}\binom {n} {2i+1}(3X_1^2)^{(n-1)/2-i}(-2Y_1^2)^i=L_n(\alpha,\beta).
\end{equation}
Therefore, by \eqref{eq.2.15} and \eqref{eq.2.19}, we get
\begin{equation}\label{eq.2.20}
\left\lvert L_n(\alpha,\beta)\right\rvert=p^{r-s}.
\end{equation}

If $s>0$, by \eqref{eq.2.14}, \eqref{eq.2.18} and \eqref{eq.2.20}, then the Lehmer number $L_n(\alpha,\beta)$ has no primitive divisors. Therefore, since $n$ is an odd prime, by Lemmas \ref{eq.2.3} and \ref{eq.2.4}, we find from \eqref{eq.2.18} that $n\in\{3,5\}.$

When $n=3$, by \eqref{eq.2.14} and \eqref{eq.2.15}, we have
\begin{equation}\label{eq.2.21}
9X_1^2-2p^{2s}=\pm p^{r-s}.	
\end{equation}
Notice that $p>3$, $s>0$ and $\gcd(X_1,Y_1)=\gcd(X_1,p^s)=1$. We see from \eqref{eq.2.21} that $r-s=0$ and
\begin{equation}\label{eq.2.22}
9X_1^2-2p^{2s}=\pm 1.
\end{equation}
Further, since $2\nmid X_1$ and $9X_1^2-2p^{2s}\equiv 1-2\equiv -1 \pmod{8}$, by \eqref{eq.2.22}, we get
\begin{equation}\label{eq.2.23}
9X_1^2-2p^{2s}=-1.
\end{equation}
But, since $(2/3)=-1$, where $(*/*)$ is the Legendre symbol, \eqref{eq.2.23} is false. So, we have no solutions for $n=3$.

When $n=5$, by \eqref{eq.2.14} and \eqref{eq.2.15}, we have
\begin{equation}\label{eq.2.24}
45X_1^4-60X_1^2p^{2s}+4p^{4s}=\pm p^{r-s}.
\end{equation}
If $r-s>0$, since $p>3$, then from \eqref{eq.2.24} we get $p=5$ and
\begin{equation*}
9X_1^4-12\cdot 5^{2s}X_1^2+4\cdot 5^{4s-1}=\pm 5^{r-s-1},
\end{equation*}
whence we obtain $r-s=1$ and
\begin{equation}\label{eq.2.25}
9X_1^4-12\cdot 5^{2s}X_1^2+4\cdot 5^{4s-1}=\pm 1.
\end{equation}
Further, since $9X_1^4\equiv 1 \pmod{4}$, the right side of \eqref{eq.2.25} is equal to 1. However, since $5\nmid X_1$ and $9X_1^4\equiv 9\equiv -1 \pmod{5},$ the right side of \eqref{eq.2.25} should be equal to -1, a contradiction. So we have $r-s=0$ and
\begin{equation}\label{eq.2.26}
45X_1^4-60X_1^2p^{2s}+4p^{4s}=\pm 1.
\end{equation}
Similarly, since $45X_1^4\equiv 1 \pmod{4}$ and $4p^{4s}\equiv -1 \pmod{5}$, \eqref{eq.2.26} is false. This implies that we have no solutions for $n=5$.

By the above analysis, we get $s=0$. Then, by \eqref{eq.2.14}, we have $Y_1=1$. Therefore, by \eqref{eq.2.2}, \eqref{eq.2.12} and \eqref{eq.2.13}, we obtain \eqref{eq.1.7} and \eqref{eq.1.8}. Thus, the theorem is proved.
\section{Proof of Theorem \ref{theo:1.2}}
For fixed $d$ with \eqref{eq.xx} and $n$ odd prime, we now assume that \eqref{maineq} has two distinct solutions $(x,y)=(x_1,y_1)$ and $(x_2,y_2)$. Then, by Theorem \ref{theo:1.1}, we have 
\begin{equation}\label{eq.3.1}
\begin{aligned}
&d=\left\lvert \sum\limits_{i=0}^{(n-1)/2}\binom {n} {2i+1}(3a^2)^{(n-1)/2-i}(-2)^i\right\rvert\\
&=\left\lvert \sum\limits_{i=0}^{(n-1)/2}\binom {n} {2i+1}(3b^2)^{(n-1)/2-i}(-2)^i\right\rvert,
\end{aligned}
\end{equation}
\begin{equation}\label{eq.3.2}
y_1=3a^2+2,\,\,y_2=3b^2+2,\,\,a,b\in\mathbb{N},\,\,2\nmid ab.
\end{equation}
Since $(x_1,y_1)\neq (x_2,y_2)$, we have $y_1\neq y_2$. Therefore, without loss of generality we may assume that $y_1<y_2$. Then, by \eqref{eq.3.2}, we get $a<b$.

Since $n$ is an odd prime, we have $n \mid \binom {n} {2i+1}$ for $i=0,\cdots,(n-3)/2$. Hence, since $n\nmid 2^{(n-1)/2}$, we see from \eqref{eq.3.1} that
\begin{equation*}
\sum\limits_{i=0}^{(n-1)/2}\binom {n} {2i+1}(3a^2)^{(n-1)/2-i}(-2)^i=\sum\limits_{i=0}^{(n-1)/2}\binom {n} {2i+1}(3b^2)^{(n-1)/2-i}(-2)^i,
\end{equation*}
whence we get
\begin{equation}\label{eq.3.3}
\sum\limits_{i=0}^{(n-3)/2}\binom {n} {2i+1}\left(\dfrac{(3b^2)^{(n-1)/2-i}-(3a^2)^{(n-1)/2-i}}{3b^2-3a^2}\right)(-2)^i=0.
\end{equation}
Let $X=3b^2$ and $Y=3a^2$. Then \eqref{eq.3.3} can be rewritten as
\begin{equation}\label{eq.3.3A}
\sum\limits_{i=0}^{(n-3)/2}\binom {n} {2i+1}\left(\dfrac{X^{(n-1)/2-i}-Y^{(n-1)/2-i}}{X-Y}\right)(-2)^i=0.
\end{equation}
By \eqref{eq.3.3A}, we have $n>3$ and
\begin{equation}\label{eq.3.4}
2\,\ \biggl |\,   \dfrac{X^{(n-1)/2}-Y^{(n-1)/2}}{X-Y}.
\end{equation}
Since $2\nmid XY$ by \eqref{eq.3.2}, we see from \eqref{eq.3.4} that $2\mid (n-1)/2$. Further let $2^y\mid\mid n-1$. Then we have $y\ge 2$ and
\begin{equation}\label{eq.3.4A}
	2^{y-1}\,\ \biggl |\biggl |\,\binom {(n-1)/2} {1}Y^{(n-3)/2}.
\end{equation}
Let $2^{r_j}\mid \mid j$ for $j>1$. Since $j\ge 2^{r_j}$, we have $r_j\le (\log j)/(\log2)\le j-1$. Since $X-Y\equiv 3a^2-3b^2\equiv 0 \pmod{2^3}$, we get
\begin{equation}\label{eq.3.4B}
\begin{aligned}
&\binom{(n-1)/2}{j}(X-Y)^{j-1}Y^{(n-1)/2-j}\\
&\equiv \left(\dfrac{n-1}{2}\right)Y^{(n-1)/2-j}\binom {(n-3)/2}{j-1}\dfrac{(X-Y)^{j-1}}{j}\\
&\equiv 0 \pmod{2^y},\,\,j>1.
\end{aligned}
\end{equation}
Hence, since
\begin{equation*}
\dfrac{X^{(n-1)/2}-Y^{(n-1)/2}}{X-Y}=\sum\limits_{j=1}^{(n-1)/2}\binom {(n-1)/2} {j}(X-Y)^{j-1}Y^{(n-1)/2-j},
\end{equation*}
we obtain from \eqref{eq.3.4A} and \eqref{eq.3.4B} that
\begin{equation}\label{eq.3.5}
2^{y-1} \biggl |\biggl | \dfrac{X^{(n-1)/2}-Y^{(n-1)/2}}{X-Y}=\dfrac{(3b^2)^{(n-1)/2}-(3a^2)^{(n-1)/2}}{3b^2-3a^2}.
\end{equation}
On the other hand, let $2^{\delta_i}\mid \mid 2i$ for $i\ge 1$. Then we have
\begin{equation}\label{eq.3.6}
\delta_i\le \dfrac{\log(2i)}{\log2}\le i,\,\,i\ge 1.
\end{equation}
By \eqref{eq.3.6}, we get
\begin{equation}\label{eq.3.7}
\binom {n} {2i+1}(-2)^i\equiv n(n-1)\binom {n-2} {2i-1}\dfrac{(-2)^i}{2i(2i+1)}\equiv 0 \pmod{2^y}, \,\,i\ge 1.
\end{equation}
Therefore, since $2\nmid n$, we find from \eqref{eq.3.5} and \eqref{eq.3.7} that \eqref{eq.3.3} is false. It implies that, under the assumption of Theorem \ref{theo:1.1}, \eqref{maineq} has at most one solution $(x,y)$. The theorem is proved.

\section{Proof of Theorem \ref{theo:1.3}}

\begin{lemma}[\cite{Leh}]\label{lem.4.1}
If $n$ is an odd prime and $r$ is a prime divisor of the Lehmer number $L_n(\alpha,\beta)$, then $r\equiv \pm 1 \pmod{2n}$.
\end{lemma}
\textbf{Proof of Theorem \ref{theo:1.3}} By Lemma \ref{lem.2.2}, if $(x,y)$ is a solution of \eqref{maineq}, then $x,y$ and $d$ satisfy \eqref{eq.2.1} and \eqref{eq.2.2}. Let $\alpha,\beta$ be defined as in \eqref{eq.2.16}. Then $(\alpha,\beta)$ is a Lehmer pair with the parameters \eqref{eq.2.18}. Further, let $L_m(\alpha,\beta)$ $(m=1,2,\cdots)$ be the corresponding Lehmer numbers. By \eqref{eq.2.13} and \eqref{eq.2.19}, we have
\begin{equation}\label{eq.4.1}
d=Y_1 |L_n(\alpha,\beta)|.
\end{equation}

Since $n$ is an odd prime and every odd prime divisor $q$ of $d$ satisfies $q\not\equiv \pm 1 \pmod{n}$, by Lemma \ref{lem.4.1}, we get from \eqref{eq.4.1} that
\begin{equation}\label{eq.4.2}
|L_n(\alpha,\beta)|=1,
\end{equation}
and
\begin{equation}\label{eq.4.3}
Y_1=d.
\end{equation}
We see from \eqref{eq.4.2} that the Lehmer number $L_n(\alpha,\beta)$ has no primitive divisors. Therefore, using the same method as in the proof of Theorem \ref{theo:1.1}, by Lemmas \ref{lem.2.3} and \ref{lem.2.4}, we can deduce from \eqref{eq.4.2} that $n\in\{3,5\}$.

When $n=3$, by \eqref{eq.2.19}, \eqref{eq.4.2} and \eqref{eq.4.3}, we have
\begin{equation}\label{eq.4.4}
9X_1^2-2d^2=\pm 1.
\end{equation}
Since $n=3$ and every odd prime divisor $q$ of $d$ satisfies $q\not\equiv \pm 1 \pmod{3}$, $q$ can only be equal to 3. However, by \eqref{eq.4.4}, it is impossible. Hence, $d$ must be a power of 2. Then \eqref{eq.4.4} reduces to the equation
\begin{equation}\label{eq.4.4.1}
X^2+1=2^{2k+1},\,\,X=3X_1,\,\,k\ge 0, 
\end{equation}
or
\begin{equation}\label{eq.4.4.2}
X^2-1=2^{2k+1},\,\,X=3X_1,\,\,k\ge 0. 
\end{equation}
By \cite{Leb}, we see that \eqref{eq.4.4.1} has no solution. Since $\gcd(X+1,X-1)=2$, we get from \eqref{eq.4.4.2} that $X-1=2$ and $k=1$. It follows that the equation has only the solution $(X,k)=(3,1).$
Therefore, it is easy to get $X_1=1$ and $d=2$. Thus, \eqref{maineq} has only the solution $(x,y,d,n)=(21,11,2,3)$ in this case.

When $n=5$, by \eqref{eq.2.19}, \eqref{eq.4.2} and \eqref{eq.4.3}, we have
\begin{equation}\label{eq.4.5}
45X_1^4-60X_1^2d^2+4d^4=\pm 1.
\end{equation}
But, since $2\nmid X_1$, $45X_1^4\equiv 1 \pmod{4}$, $5\nmid d$ and $4d^4\equiv -1\pmod{5}$, \eqref{eq.4.5} is false. The theorem is proved.

\section{Proof of Theorem \ref{theo:1.4}}
For any algebraic number $\theta$ of degree $\ell$ over $\mathbb{Q}$, let $h(\theta)$ be the absolute logarithmic height of $\theta$ by the formula
\begin{equation*}
h(\theta)=\frac{1}{\ell} \left(\log|a|+ \sum_{j=1}^{\ell} \log\max\big\{1,|\theta^{(j)}|\big\}\right)
\end{equation*} 
where $a$ is the leading coefficient of the minimal polynomial of $\theta$ over $\mathbb{Z}$ and $\theta^{(j)}$ $(j=1,\cdots,\ell)$ are all the conjugates of $\theta$. Further, let $\log\theta$ be any determination of the logarithm of $\theta$.
\begin{lemma}[Appendix of \cite{BHV}]\label{lem.5.1}
Let $\theta$ be a complex algebraic number with $|\theta|=1,$ and $\theta$ is not a root of unity. Let $b_1$, $b_2$ be positive integers, and let $\Lambda=b_1\log\theta-b_2\pi\sqrt{-1}.$ Then we have
\begin{equation*}
\log|\Lambda|>-(9.03H^2+0.23)(Dh(\theta)+25.84)-2H-2\log H-0.7D+2.07,
\end{equation*} 
where $D=[\mathbb{Q}(\theta):\mathbb{Q}]/2,$ $H=D(\log B-0.96)+4.49$, $B=\max\{13,b_1,b_2\}$.
\end{lemma}
\textbf{Proof of Theorem \ref{theo:1.4}} By Lemma \ref{lem.2.2}, if $(x,y)$ is a solution of \eqref{maineq}, then
\begin{equation}\label{eq.5.1}
d=\dfrac{1}{2\sqrt{2}}|\alpha^n-\beta^n|,
\end{equation}
where $\alpha,\beta$ are defined as in \eqref{eq.2.16}. By \eqref{eq.2.2} and \eqref{eq.2.16}, we have
\begin{equation}\label{eq.5.2}
|\alpha|=|\beta|=\sqrt{y}.
\end{equation}
Let $\theta=\alpha/\beta$. By \eqref{eq.5.2} and \eqref{eq.2.17-A}, is a complex algebraic number with $|\theta|=1$, $\theta$ is not a root of unity and
\begin{equation}\label{eq.5.3}
h(\theta)=\dfrac{1}{2}\log y.
\end{equation}

By \eqref{eq.5.1} and \eqref{eq.5.2}, we have
\begin{equation}\label{eq.5.4}
d=\dfrac{1}{2\sqrt{2}}|\beta^n|\left\lvert\left(\dfrac{\alpha}{\beta}\right)^n-1\right\lvert=\dfrac{1}{2\sqrt{2}}y^{n/2}|\theta^n-1|.
\end{equation}
It is well known that, for any complex number $z$, we have either $|e^z-1|\ge \dfrac{1}{2}$ or $|e^z-1|\ge \dfrac{2}{\pi}|z-t\pi\sqrt{-1}|$ for some integers $t$ (see \cite{St}). Put $z=n\log \theta$. We get either
\begin{equation}\label{eq.5.5}
|\theta^n-1|\ge \dfrac{1}{2},
\end{equation}
or
\begin{equation}\label{eq.5.6}
|\theta^n-1|\ge \dfrac{2}{\pi}|n\log \theta-t\pi\sqrt{-1}|,\,\,t\in\mathbb{N},\,\,t\le n.
\end{equation}
If \eqref{eq.5.5} holds, since $d>8\sqrt{2}$, then from \eqref{eq.5.4} we obtain $y^n\le 32d^2<2^{3/2}d^3$ and the theorem is true. So we just have to worry about the case \eqref{eq.5.6}.

Let
\begin{equation}\label{eq.5.7}
\Lambda=n\log \theta-t\pi\sqrt{-1}.
\end{equation} 
By \eqref{eq.5.4}, \eqref{eq.5.6} and \eqref{eq.5.7}, we have
\begin{equation}\label{eq.5.8}
d\ge \dfrac{y^{n/2}}{\pi\sqrt{2}}|\Lambda|.
\end{equation}
If $y^n\ge 2^{3/2}d^3$, then from \eqref{eq.5.8} we get
\begin{equation*}
\pi \ge y^{n/6}|\Lambda|,
\end{equation*}
whence we obtain
\begin{equation}\label{eq.5.9}
\log \pi \ge \dfrac{n}{6}\log y+\log |\Lambda|.
\end{equation}

Notice that $[\mathbb{Q}(\theta):\mathbb{Q}]=2$, $n\ge t$ and $n>228000$. Applying Lemma \ref{lem.5.1} to \eqref{eq.5.7}, by \eqref{eq.5.3}, we have
\begin{equation}\label{eq.5.10}
\log |\Lambda|>-(9.03H^2+0.23)(\dfrac{1}{2}\log y+25.84)-2H-2\log H+1.37,
\end{equation}
where
\begin{equation}\label{eq.5.11}
H=\log n+3.53.
\end{equation}
The combination of \eqref{eq.5.9} and \eqref{eq.5.10} yields
\begin{equation}\label{eq.5.12}
(9.03H^2+0.23)\left(0.5+\dfrac{25.84}{\log y}\right)+\dfrac{2H+2\log H}{\log y}>\dfrac{n}{6}.
\end{equation}
Further, by \eqref{eq.2.2}, we have $y\ge 5$. Hence, by \eqref{eq.5.11} and \eqref{eq.5.12}, we get
\begin{equation}\label{eq.5.13}
\begin{aligned}
& 99.36(9.03(\log n+3.53)^2+0.23)+7.50(\log n+3.53\\
&+\log(\log n+3.53))=99.36(9.03H^2+0.23)\\
&+7.50(H+\log H)>n.
\end{aligned}
\end{equation} 
However, by \eqref{eq.5.13}, we calculate that $n<228000$, a contradiction. Thus, if $n>228000$ and $d>8\sqrt{2}$, then $y^n<2^{3/2}d^3.$ The theorem is proved.

\subsection*{Acknowledgements}
This paper was partially written when the second author participated in the workshop titled \textquotedblleft Effective Methods for Diophantine Problems\textquotedblright \, on 18-22 June 2018 in Lorentz Workshop Center (Leiden University), Netherlands. He would like to thank organizers Professors Attila B\'{e}rczes, Bas Edixhoven, Kalman Gy\H{o}ry and Robin Je Dong for giving an opportunity to participate in this excellently organized workshop, giving accommodation support and their kind hospitality. The authors would like to Dr. Paul Voutier for useful discussions and to Dr. Angelos Koutsianas for sharing the final version of his paper.

\end{document}